\documentclass[12pt]{article}
\usepackage{amssymb}
\usepackage{amsmath}
\newtheorem{theo}{Theorem}

\newtheorem{prop}[theo]{Proposition}
\newtheorem{lemm}[theo]{Lemma}
\newtheorem{coro}[theo]{Corollary}

\newtheorem{conj}[theo]{Conjecture}
\newtheorem{quest}[theo]{Question}

\def \P{{\mathbb P}}

\begin{document}

\title{On families of lagrangian tori on hyperk\"ahler manifolds}

\author{Ekaterina Amerik, Fr\'ed\'eric Campana}

\date{}

\maketitle

\section{Introduction}

Recall that an irreducible holomorphic symplectic manifold is a
simply-connected compact K\"ahler manifold $X$ such that the space
of holomorphic $2$-forms $H^{2,0}(X)$ is generated by a symplectic,
that is, nowhere
degenerate, form $\sigma$. It is well-known that such manifolds are
exactly the compact hyperk\"ahler manifolds from differential geometry.
There are two series of known examples, namely the Hilbert scheme 
$Hilb^{[n]}(X)$ of length $n$ subschemes of a $K3$ surface $X$ and the
Kummer variety associated to an abelian surface $A$ (recall that this 
is defined as the fiber over zero of the summation map $\Sigma: Hilb^{[n]}(A)\to A$), and two
sporadic examples in dimensions 6 and 10 due to O'Grady. All known
irreducible holomorphic symplectic manifolds are deformations of those.

In order to study the classification problem for irreducible holomorphic symplectic 
manifolds and their geometry,
it is important to understand how such an $X$ can fiber over lower-dimensional varieties. 
Several results are known in this direction.
The first one, very striking, has been obtained by D. Matsushita:

\begin{theo}\label{mats-fibr} (\cite{M1},\cite{M2}) Let $X$ be an 
irreducible
holomorphic symplectic variety and let $f:X\to B$ be a map with connected
fibers from $X$ onto a normal complex space $B$, $0<dim(B)<n$. 
Then all fibers of
$f$ are lagrangian (in particular, they are of dimension $n$), 
and the smooth fibers are tori\footnote{More precisely, this is proved for
$X$ and $B$ projective in \cite{M1} and it is remarked in \cite{M2} that the arguments work in the 
K\"ahler case (that is, with $X$ an arbitrary irreducible holomorphic symplectic manifold
and $B$ a normal K\"ahler space). In general,
by \cite{V2} $B$ is birational to a K\"ahler manifold $B'$, and it follows from the 
Kodaira embedding
theorem that $B'$ is projective. Replacing
$f:X\to B$ by $f':X\dasharrow B'$ and considering an ample $H'$ on $B'$, one 
shows (as in \cite{AC}, Prop. 3.1)
that $(f')^*H'$ is isotropic with respect to the 
Beauville-Bogomolov form. One then proves
as in \cite{M1} that a general fiber of $f'$, and therefore of $f$, is 
lagrangian. By \cite{M1bis}, $f$ is equidimensional, and therefore $B$ is normal 
K\"ahler by \cite{V1}.}. 
\end{theo}

Strong restrictions on $B$ are known. First, the results of Varouchas (\cite{V1},\cite{V2}) 
imply the existence of a smooth K\"ahler model $B'$ of $B$, and $B'$ does 
not carry any nonzero holomorphic $2$-form (since its lift to $X$ would not be proportional to $\sigma$), so
Kodaira's embedding theorem shows that $B'$ is projective. 

Much more is true. Matsushita proved that the base $B$ of such a fibration is
``very similar'' to $\P^n$. Later, Hwang has shown in \cite{H} that 
if $B$ is smooth, then indeed $B\cong\P^n$ (the general case is open).

\medskip

In \cite{AC} we have partially extended this to the meromorphic setting.
Since there are always meromorphic fibrations by complete intersections
of hypersurfaces if $X$ is projective, we assumed that the general 
fiber of 
$f:X\dasharrow B$ is not of general type; in this case, we have proved, 
among other things,
that $dim(B)\leq n$ (proposition 3.4 of \cite{AC}).

When $X$ is
non-projective, our assumption on the fibers is automatically satisfied
since the total space of a fibration with Moishezon base and general
type fibers is Moishezon, see \cite{U}, 2.10.

In fact the general fiber of $f$ is not of general type if and only if
$f$ is defined by sections of a line bundle with zero Beauville-Bogomolov
square (\cite{AC}, Proposition 3.1). This is related to a famous and difficult conjecture 
about lagrangian fibrations:

\begin{conj} (``Lagrangian conjecture''): Let $L$ be a nontrivial nef line bundle on $X$,
such that $q(L)=0$. Then some power of $L$ is base-point-free, that is, the sections of
some power of $L$ define a holomorphic lagrangian fibration.
\end{conj}

Under the additional assumption that $X$ is covered by curves $C$ such that  
$LC=0$, this has been proved by Matsushita \cite{M3}. The main idea of the proof is to use the 
{\it nef reduction} from \cite{BCE...}, which in
this case yields a non-trivial meromorphic fibration.

Recently, Beauville has asked the following question:

\begin{quest}\label{bq} Let $X$ be an irreducible holomorphic symplectic variety
and $A\subset X$ a lagrangian torus. Is it true that $A$ is a fiber
of a (meromorphic) lagrangian fibration?
\end{quest}

This is plausible since it is known that the deformations of $A$ in $X$ are 
unobstructed and that
the smooth ones are again lagrangian tori. The symplectic
form defines an isomorphism between the normal
bundle of a lagrangian torus and its cotangent bundle, so this bundle is
trivial. Hence, locally in a neighbourhood of $A$, there is a lagrangian 
fibration, and the question is whether it globalizes. Another reformulation
is as follows: one knows that through a general point of $X$ passes a finite
number $d$ of deformations of $A$, and one wants to know whether $d=1$.

When $X$ is non-projective, the affirmative answer has been given
by Greb,  Lehn and Rollenske in \cite{GLR}. Moreover, they have proved that 
if the pair $(X, A)$ can be deformed to a pair $(X', A')$ where $X'$ 
is non-projective 
and $A'$ is still a lagrangian torus on $X'$, then the
answer is also positive (because the lagrangian fibration on $X'$ deforms 
back to $X$).
Finally, they have observed that by deformation theory of hyperk\"ahler 
manifolds, 
the existence of a non-projective deformation of the pair $(X, A)$ is 
equivalent to the 
existence of an effective divisor $D$ on 
$X$ such that its restriction to $A$ is zero.

In \cite{A}, the first author has given a very simple
solution of Beauville's problem 
in dimension 4, by combining a lemma from linear
algebra with the above-mentioned proposition from \cite{AC}.

In \cite{HW}, the existence of an effective divisor $D$ restricting
trivially to $A$ is established, 
so that the affirmative answer to question \ref{bq} is obtained in
full generality.
The argument, which proceeds by
the study of the monodromy action on the total space of the family of 
deformations of $A$, uses some highly non-trivial finite group theory.

This solution of Beauville's question, as well as
several other important proofs in holomorphic symplectic geometry,
depends in an essential way on deforming from projective to non-projective 
data.

The purpose of this note is to give another, very simple, proof in the
non-projective case (the original one from \cite{GLR} uses results about 
algebraic reduction of hyperk\"ahler
manifolds from \cite{COP} and involves some case-by-case analysis) and to 
explain
a possible algebro-geometric
approach to Beauville's question (not relying on deformations or on group
theory). We aim to explicitely construct a linear system which, 
possibly after some flops, would give us a Lagrangian fibration 
$f:X\to B$ with generic fiber meeting the generic member or our Lagrangian 
family $A$ along a
positive-dimensional subvariety. Then hopefully one can show 
that this implies that $A$ is a fiber of $f$ 
\footnote{For instance, let $f:X\to B$ and $g:X\to C$ be two Lagrangian
fibrations. Then either $f=g$, or the generic fibers of $f$ and $g$ have finite
intersections. This is seen by using the fact that otherwise 
the Beauville-Bogomolov 
form
vanishes on the subspace generated by divisors coming from $B$ and $C$.}.

Unfortunately our algebro-geometric argument does not, at the moment,
give a general answer: we can produce a lagrangian fibration in a special case 
(corollary \ref{bymmp}), and only a weaker statement (theorem \ref{zerointersection}) is 
obtained in general.

\section{The non-projective case}

Let us start with the following observation, which is an immediate 
consequence of the argument in \cite{C2}, Prop. 2.1. 

\begin{prop}\label{projectivity} A smooth lagrangian subvariety of an irreducible 
holomorphic symplectic manifold is projective.
\end{prop}

By \cite{C1}, 
there is an almost holomorphic fibration $\phi_S$ associated to any compact 
and covering 
family $Z_s,\ s\in S$, of subvarieties in $X$: the fiber of $\phi_S$
through a sufficiently general point $x\in X$ consists of all points
$y$ such that there is a chain of subvarieties from $S$ joining
$x$ and $y$. Moreover, again by \cite{C1}, if $Z_s$ are projective varieties, 
then so are
the fibers of $\phi_S$. 

Throughout the paper, we shall also need the following simple lemma, which may be seen 
as an analogue of our problem for tori:

\begin{lemm}\label{fibtori} Let $Z$ be an irreducible subvariety 
of a complex torus $T$. Assume that an open neighborhood $U$ of $Z$ 
admits a proper holomorphic fibration $g_U:U\to V$ having $Z$ as one of 
its fibers. Then $Z$ is a complex subtorus of $T$, so that the local fibration yields
a global fibration $g: T\to B$. 
\end{lemm}

{\it Proof:} The generic fiber of $g_U$ is a smooth connected complex
 submanifold $Z'$ of $T$. By adjunction, its canonical bundle is trivial. 
By Ueno's theorem (\cite{U}, theorem 10.3), $Z'$ is a complex torus, thus 
so is $Z$. Indeed, $Z$ is a translate of $Z'$, since the translates
 of $Z'$ form a connected component of the Barlet-Chow space of $T$.

\begin{theo}\label{nonproj}
Let $X$ be a non-algebraic irreducible holomorphic symplectic manifold
and $A$ a lagrangian torus on $X$. Then $X$ admits an almost holomorphic
lagrangian fibration such that $A$ is a fiber.
\end{theo}

{\it Proof:} Consider the family of deformations $A_t,\ t\in T$ of $A$
and the fibration $\phi_T$ associated to $A_t$ (note that 
$T$ is compact since $X$ is compact K\"ahler). We claim that its relative
dimension is $n$, that is, $A_t$ is
a fiber of $\phi_T$ for general $t$. Indeed, the relative dimension is
obviously at least $n$. If it is $2n$, that is, $\phi_T$ is a constant 
map, then $X$, being a fiber of $\phi_T$, is projective because the $A_t$
are projective. If it is strictly between $n$ and $2n$, consider, for
a general $x$, all tori $A_1,\dots, A_d$ passing through $x$. By 
assumption, $d\geq 2$. The tori $A_1,\dots, A_d$ are contained in the fiber $F_x$ of $\phi_T$
through $x$, in particular, the dimension of the subspace $V_x\subset T_{X,x}$ generated by 
$T_{A_1, x},\dots,T_{A_d, x}$ is strictly less than $2n$.

From the fact that in a neighbourhood of a general torus, our family is a fibration,
we easily deduce that if $x\in X$ is general, the tori $A_1,\dots, A_d$ are not tangent
to each other at $x$, meaning that their intersection has only one component $Z_x$ 
through $x$ and $n\geq dim(Z_x)=dim T_{Z_x, x}= dim \cap_iT_{A_i, x}$. But $\cap_iT_{A_i, x}$
is exactly the $\sigma$-orthogonal to $V_x$ and thus it is strictly positive-dimensional, meaning
that so is the component of the intersection of $A_i$ through $x$. One thus obtains
a meromorphic fibration of $X$ by such components; call them $E_x$. If one knows that
$E_x$ is not of general type, then by Proposition 3.4 of \cite{AC} one concludes that
$dim(E_x)=n$. By definition of $E_x$, this means that $d=1$ and the family $A_t$ fibers $X$, q.e.d.. 

The fact that $E_x$ is not of general type
is easily deduced by induction on $d$: indeed we know that the family $A_t$ gives a local fibration
near its general member, so the same must be true for intersections $A_t\cap A_s$ for fixed general $A_s$ and
varying general $A_t$ intersecting $A_s$. But by lemma \ref{fibtori}, a torus can only be locally fibered in subtori. If $x$ is general, then so are
$A_1$ and $A_2$, so the component of $A_1\cap A_2$ through $x$ is a torus. Continuing in this way,
we conclude that also $A_1\cap\dots \cap A_d$ is a torus. This finishes the proof.

\section{The projective case: some results}

Consider the family of lagrangian tori $A_t, t\in T$, which are deformations of a certain lagrangian torus
$A\subset X$. Recall that these deformations cover $X$, and that there exists a number 
$d$ such that exactly $d$ members of this family pass through a generic point of $X$. 
Assuming that this does not yield a meromorphic fibration, that is, $d>1$, 
we are going to construct two large families of subvarieties of $X$ of complementary dimensions $e$ and $2n-e$,
such that the intersection number of the corresponding cycles in the cohomology is zero.

Recall that the intersection $A_s\cap A_t$ has no zero-dimensional components (this follows from
the fact that our family is a local fibration, see \cite{A}, lemma 1). For fixed $t$, the tori 
intersecting $A_t$ form a finite number of irreducible families $S_1,\dots, S_N\subset T$. Fix one of them 
and call it $S$: for $s\in S$ general, $A_s\cap A_t$ is an equidimensional (\cite{A}, lemma 1) union of 
disjoint subtori
by lemma \ref{fibtori}. We shall see very soon that the all these subtori are translates of each other in $A_t$. 

Set $e=e_S$ the dimension of $A_s\cap A_t$ for general $s\in S$. Denote by $E_{s,t}$ a 
component of $A_s\cap A_t$, and by $Z_{s,t}$ the union
of $A_s$, $s\in S$, which is thus of dimension $2n-e$.

Our main result is the following theorem.

\begin{theo}\label{zerointersection}
 $[E_{s,t}]\cdot[Z_{s,t}]=0$ in the cohomologies of $X$.
\end{theo}

Before proving the theorem, let us state a corollary which partially answers Beauville's question
in the case $e=1$ (i.e. in the case when there exists an $S$ as above such that $e_S=1$). Notice
that the ``opposite'' case when $e=n-1$ for all $S$ has a completely elementary treatment
(\cite{A}, remark 4). Unfortunately we did not succeed in proving an analogue of this corollary for
arbitrary $e$.

\begin{coro} \label{bymmp} Suppose $e=1$. Then for some $m\geq 0$, the linear system $|mZ|=|mZ_{s,t}|$ 
gives an almost 
holomorphic map $\phi: X\dasharrow B$, which has
a holomorphic model $\phi': X'\to B'$ (and therefore is a lagrangian fibration).

\end{coro}

{\it Proof:} This would follow at once from theorem 
\ref{zerointersection} and the main theorem of \cite{M3}
if the divisor $Z$ had been nef: indeed $X$ is covered by curves $C=E_{s,t}$
with $CZ=0$. In general, one cannot affirm that $Z$ is nef. But it is
mobile, that is, has no base components. Such a divisor on a 
hyperk\"ahler manifold can be made nef by a sequence of flops
(see \cite{M4}). That is, there is a birational transformation
$h: X\dasharrow X'$ which is an isomorphism in codimension one, and
such that $h_*Z$ is nef.

To see that the sections of some power of $h_*Z$ give a lagrangian
fibration  $\phi': X'\to B'$, it suffices now to remark that a
general member of the family of deformations of $E_{s,t}$ does
not intersect the locus $W\subset X$ where $h$ is not an isomorphism.
This follows easily from the observation that a neighbourhood 
$U_s\subset X$ of $A_s$ is fibered not only by lagrangian tori which
are small deformations of $A_s$, but 
also by deformations of the subtorus $E_{s,t}$: the local lagrangian
fibration 
$f_s: U_s\to B_s$ factors through $g_s: U_s\to D_s$ which has
$E_{s,t}$ for a fiber. The locus $W$, being of codimension at least
two, cannot dominate $U_s$ and so we can deform $E_{s,t}$ away
from it and then take the image by $h$ to obtain a dominating family of 
curves which do not intersect $h_*Z$. With all this done, we conclude
by \cite{M3}.

\medskip

The proof of the theorem uses several preliminary lemmas related to
the following construction (``the characteristic foliation'').

Consider the local lagrangian fibration $f_s: U_s\to B_s$ in
a neighbourhood of $A_s$. Then $A_t\cap U_s$ projects onto 
a codimension $e$ analytic subvariety $C_s\subset U_s$. Take 
$D^0=f_s^{-1}(C_s)$. This is the ``principal branch" of the subvariety
$Z_{s,t}$ in the neighbourhood of $A_s$. It is clear that we can
move $E_{s,t}$ away from $D^0$ by replacing it on some 
neighbouring $A_{s'}$ with $s'\not\in C_s$. But there are other
branches $D^1, \dots, D^l$ of $Z_{s,t}$ intersecting $A_s$ and all of its
neighbours, and we must show that $E_{s,t}$ can be moved away from them, too.

For this, we look at the kernel of the restriction of the symplectic form 
$\sigma$ to the smooth part of $Z_{s,t}$: this is a distribution 
$\cal F$ of rank $e$, often called characteristic foliation in the
literature.

\begin{lemm} \label{trivial} The restriction of $\cal F$ to $A_s$ is trivial
(as a subbundle of $T_{A_s}$).

\end{lemm}

{\it Proof:} The symplectic form $\sigma$ defines in the usual
way an isomorphism between the restriction of $\cal F$ to $D_0$ and
the conormal bundle to $D_0$ in $X$: indeed the last one is the
kernel of the restriction map from the dual of $T_X|_{D_0}$ to the dual
of $T_{D_0}$. We claim that the restriction of the normal bundle 
$N_{D_0, X}$ to $A_s$ is trivial.
This is clear from the exact sequence
$$0\to N_{A_s, D_0}\to N_{A_s, X}\to N_{D_0, X}|_{A_s}\to 0:$$
the first term is a trivial vector bundle since $D_0$ is a fibration, the 
second one is
trivial because $A_s$ is lagrangian and therefore $N_{A_s, X}$ is isomorphic
to the cotangent bundle of $A_s$, and the third one is trivial because
it therefore has $e$ everywhere linearly independent global sections.

\begin{coro} The distribution ${\cal F}$ is tangent to a fibration of
$Z_{s,t}$ by deformations of $E_{s,t}$.
\end{coro}

{\it Proof:} The restriction of ${\cal F}$ to $A_s$ has at least one algebraic
leaf, that is, $E_{s,t}$ (being a component of the intersection of lagrangian
$A_s$ and $A_t$, it is orthogonal to both). Since it is trivial as a vector bundle,
all leaves are translates of $E_{s,t}$ in  $A_s$. Since $A_s$ is generic, all leaves
on $Z_{s,t}$ are algebraic.

\medskip

The following lemma, which seems to be an explicit geometric analogue of 
Hwang--Weiss'
``pairwise integrability'' (though we have obtained it independently around the same time), is crucial.

\begin{lemm}\label{bicouvert} There is only a finite number of $Z_{s,t}$ passing through a 
general point $x\in X$. In particular, if $A_{t'}$ is a small deformation of $A_t$ which
still intersects $A_s$, then $Z_{s,t}=Z_{s,t'}$ and $A_{t'}\subset Z_{s,t}$. Therefore,
through a general translate of $E_{s,t}$ in $A_s$ (and not only in $A_t$), there is a 
lagrangian torus contained in $Z_{s,t}$.
\end{lemm}   

{\it Proof:} Indeed, each $Z_{s,t}$ is a union of tori $A_s$ and their degenerations
(which do not cover $X$).
Therefore a $Z_{s,t}$  passing through a general $x$ should contain such a torus $A_s$ 
through $x$,
and its tangent space at $x$ must be $\sigma$-orthogonal to that of a subtorus of 
dimension $e$ 
in $A_s$. But there is only a countable number of possibilities for those. Since we are
dealing with bounded families, we conclude that there are in fact only finitely many 
possibilities for $T_xZ_{s,t}$. Since $x$ is general, an application of standard results
(either the unicity theorem for solutions of differential equations on a 
suitable covering of $X$ or Sard's lemma on a suitable fibered product) shows that there is also 
only a finite number of $Z_{s,t}$ through $x$.

The second assertion follows from the first since $Z_{s,t}$ through $x$ does not deform. 
The third one follows from the second by taking closure.

\medskip

As an immediate corollary we obtain an assertion already announced before the statement of theorem \ref{zerointersection}:

\begin{coro}\label{translate} The intersection $A_s\cap A_t$ for $s$ general is a union of 
translates of
$E_{s,t}$ both in $A_s$ and in $A_t$.
\end{coro}

Indeed, this intersection must be tangent to the kernel of the restriction of $\sigma$,
and as we have seen above, this kernel is tangent to the fibration given by the translates. 

\medskip

{\it Proof of theorem \ref{zerointersection}:} The subvariety $Z_{s,t}$ is the union of 
$A_s$
where $A_s$ varies in an irreducible ($n-e$)-parametric family $S$ of tori intersecting 
$A_t$. As we have just seen,
$A_s\cap A_t$ for $s$ general is a union of 
translates of
$E_{s,t}$ in $A_t$: our first claim is that the same is true for for any $s$ (but the union
can a priori be infinite, that is, of greater dimension than $e$). This is easily seen
from the factorisation of the local fibration $f_s: U_s\to B_s$ through $g_s: U_s\to D_s$
from the proof of corollary \ref{bymmp}.

Next, let $s$ be general and $u\in S$ arbitrary. We claim that $A_u\cap A_s$ is again a union
of translates of $E_{s,t}$ in $A_s$. Suppose first that  $A_u$ passes through a point 
$y\in E_{s,t}\subset A_s\cap A_t$. Then by what we have just observed, $A_u$ contains 
$E_{s,t}$
(indeed $A_u\cap A_t$ is a union of translates of $E_{s,t}$, so if this contains a point
$y\in E_{s,t}$, then it contains the whole of $E_{s,t}$). Now we can repeat the same
 argument
supposing that $A_u$ passes through a point 
$y'$ on a translate $E'$ of $E_{s,t}$ in $A_s$: indeed by lemma \ref{bicouvert} we can find an $A_{t'}$
contained in $Z_{s,t}$ such that  $E'=E_{s,t'}$ is a component of the intersection 
$A_s\cap A_{t'}$ and just replace $t$ by $t'$.
In conclusion, together with each of its points, $A_u\cap A_s$ contains the whole
translate of  $E_{s,t}$ passing through this point. This proves our second claim.

Now we are able to show that $E_{s,t}\cdot Z_{s,t}=0$. Recall that we have denoted by
$D^0$ the principal branch of $Z_{s,t}$ around $A_s$, so that $D^0$ is a union of fibers of
the local fibration $f_s$, and by $D^1,\dots, D^l$ the other branches.
It is clear that a general deformation of $E_{s,t}$ is disjoint from $D^0$: it suffices to 
move $E_{s,t}$ to a neighbouring torus not contained in $D^0$. But now we see that it is
disjoint from $D^1,\dots, D^l$ as well: indeed $D^i\cap A_s$ is a union of translates 
of $E_{s,t}$ and so replacing $E_{s,t}$ by a suitable translate in $A_s$ we can make it
disjoint from $D^i$.

This finishes the proof of theorem \ref{zerointersection}.

\medskip

As a final remark, let us notice that our arguments also prove the following
statement, which can be of independent interest:

\begin{prop}\label{cst} Let $U$ resp. $V$ be connected complex manifolds of dimension
$2n$ resp. $n$, and let $\sigma$ be a holomorphic symplectic form on $U$.
Let $\varphi:U\to V$ be a proper lagrangian fibration with fibers $A_v, v\in V$
(which are then automatically complex tori). Assume that there exists a connected 
open subset $W\subset U$ and a holomorphic map $\eta:W\to M$ with $M$ a connected 
complex manifold of dimension $n$, such that $W_m:=\eta^{-1}(m),m\in M$ is a 
connected closed Lagrangian submanifold of $U$, that the restriction 
$\varphi_m:W_m\to V$ of $\varphi$ to $W_m$ is proper for every $m$, and that
the restriction $\eta_v:W\cap A_v\to M$ is proper for each $v\in V$. 

Then the image $V_m=\varphi(W_m)$ (which is a closed irreducible complex analytic
 subset of $V$) is independent of $m\in M$.
 
 More precisely, $\varphi=q\circ \alpha$ with $q:U\to N$ and $\alpha:N\to V$ proper submersions: the fibers of $q$ are connected components of $A_v\cap W_m$ and their translates
in $A_v$. This factorisation induces a smooth foliation $F$ on $V$ such that at each point $v\in V$, $F_v$ is the projection by $\varphi_*$ of the $\sigma$-orthogonal to the tangent space to the fibre of $q$ at any point of $A_v$, and $V_m$ is a leaf of $F$ for any $m\in M$.

\end{prop}

{\it Proof:} For $v,m$ generic in $V\times M$, $A_v\cap W_m$ is smooth and (by lemma \ref{fibtori}) 
is a finite union of translated subtori $E_{v,m}$ of $A_v$. The tangent space to the subvariety 
$Z_{v,m}$ (defined in the same 
way as before) at a generic 
point $x$ of a component of $A_v\cap W_m$ is generated by $T_{A_v,x}$ and $T_{W_m,x}$, and is thus the
 $\sigma$-orthogonal (by the Lagrangian property of $A_v$ and $W_m$) of $T_{E_{v,m}, x}$. By the rigidity 
(up to translation) of subtori of $A_v$, this $E_{v,m}$ is, up to translation, independent of $m$, 
and $T_{Z_{v,m}, x'}$ is thus, for every $x'\in U\cap A_v$, the $\sigma$-orthogonal to $TE_{x',v}$, 
where $E_{x',v}$ is the translate of $E_{v,m}$ through $x'$. This shows that the tangent space of $V_m$ 
at $v$ is independent of $m$, and concludes the proof of the first claim. The rest is similar.

{}

\begin{thebibliography}{}

\bibitem[A]{A} E. Amerik: A remark on a question of Beauville about lagrangian fibrations, 
Moscow Math. J. 12 (2012), no. 4, 701-704.
\bibitem[AC]{AC} E. Amerik, F. Campana: Fibrations m\'eromorphes sur certaines vari\'et\'es \`a fibr\'e canonique trivial, Pure Appl. Math. Q. 4 (2008), no. 2, part 1, 509-545.
\bibitem[BCE...]{BCE...} Th. Bauer et al: A reduction map for nef line bundles. Complex geometry (G\" ottingen, 2000), 27-36, Springer, Berlin, 2002.
\bibitem[C1]{C1} F. Campana, Cor\'eduction alg\'ebrique d'un espace analytique faiblement 
k\"ahl\'erien compact, Invent. Math. 63 (1981), no. 2, 187-223.

\bibitem[C2]{C2} F. Campana, Isotrivialit\'e de certaines familles k\" ahl\'eriennes de vari\'et\'es non 
projectives,
Math. Z. 252 (2006), no. 1, 147-156. 

\bibitem[COP]{COP} F. Campana, Th. Peternell, K. Oguiso: Non-algebraic 
hyperk\"ahler manifolds,
J. Differential Geom. 85 (2010), no. 3, 397-424


\bibitem[GLR]{GLR} D. Greb, C. Lehn, S. Rollenske: Lagrangian fibrations on hyperkähler manifolds - On a question of 
Beauville, to appear in
Annales scientifiques de l'ENS 46, no. 3 (2013).
\bibitem[H]{H} J.-M. Hwang: Base manifolds for fibrations of projective irreducible symplectic manifolds. 
Invent. Math. 174 (2008), no. 3, 625-644.
\bibitem[HW]{HW} J.-M. Hwang, R. Weiss: Webs of Lagrangian tori in
projective symplectic manifolds, to appear in Invent. Math., 
published online in May 2012.

\bibitem[M1]{M1} D. Matsushita: On fibre space structures of a projective irreducible symplectic manifold. Topology 38 (1999), no. 1, 79-83, and addendum in 40 (2001), no.2, 431-432. 

\bibitem[M1bis]{M1bis} D. Matsushita: Equidimensionality of Lagrangian fibrations on holomorphic symplectic manifolds, Math. Res. Letters 7 (2000), no. 4, 389-391. 


\bibitem[M2]{M2} D. Matsushita: Holomorphic symplectic manifolds and Lagrangian fibrations. Monodromy and differential equations (Moscow, 2001). Acta Appl. Math. 75 (2003), no. 1-3, 117-123.
\bibitem[M3]{M3} D. Matsushita: On nef reductions of projective irreducible symplectic manifolds. Math. Z. 258 (2008), no. 2, 267-270. 
\bibitem[M4]{M4} D. Matsushita: On almost holomorphic Lagrangian fibrations, arXiv:1209.1194

\bibitem[U]{U} K. Ueno: Classification theory of compact complex
analytic spaces, Lecture Notes in Math. 439 (1975).

\bibitem[V1]{V1} J. Varouchas: Stabilit\'e de la classe des vari\'et\'es 
k\"ahl\'eriennes par certains morphismes propres,
Invent. Math. 77 (1984), no. 1, 117-127.
\bibitem[V2]{V2} J. Varouchas: Sur l'image d'une vari\'et\'e 
k\"ahl\'erienne compacte, Fonctions de plusieurs variables complexes, 
V (Paris, 1979-1985), 245-259,
Lecture Notes in Math., 1188, Springer, Berlin, 1986.


\end{thebibliography}
\end{document}